\documentclass[a4paper, 11pt,russian]{article}
\usepackage{cmap}
\usepackage[T2A]{fontenc}
\usepackage[cp866]{inputenc}
\usepackage[russian]{babel}
\usepackage{amsfonts,amsmath,amssymb,theorem}
\usepackage[unicode]{hyperref} 
\usepackage{color}
\usepackage{ifthen}

\makeatletter
\def\@seccntformat#1{\S\,\csname the#1\endcsname.\ } 
\def\@biblabel#1{#1.} 
\makeatother

\date{}

\newcommand\proofr{\par{\sc Д\,о\,к\,а\,з\,а\,т\,е\,л\,ь\,с\,т\,в\,о\,.\,~}}

\def\proofend{$\blacktriangle$\vspace{0.3em}\par}

{\par\addvspace{1mm}{\sc Доказательство\hspace{1.0ex}{#1}.} }%
{\quad$\blacktriangle$\par\addvspace{1mm}}

    \newif\ifNoRemark
    \def\addtheorem#1#2#3#4{ 
    \ifthenelse{\expandafter\isundefined\csname the#2\endcsname}{\newcounter{#2}}{}
    \newenvironment{#1}[1][\global\NoRemarktrue]
     {\par\addvspace{2mm}\noindent 
       \refstepcounter{#2}{\bf #3~\csname the#2\endcsname
      \vphantom{##1}\ifNoRemark.\ \else\ (##1).\fi}\begingroup #4}%
     {\endgroup\par\addvspace{1mm}\global\NoRemarkfalse}
    \expandafter\newcommand\csname b#1\endcsname{\begin{#1}}
    \expandafter\newcommand\csname e#1\endcsname{\end{#1}}
    }

\addtheorem{pro}{pro}{Предложение}{\sl}
\addtheorem{theorem}{theorem}{Теорема}{\sl}
\addtheorem{lemma}{lemma}{Лемма}{\sl}
\addtheorem{corol}{corol}{Следствие}{\sl}
\addtheorem{sta}{sta}{Утверждение}{\sl}

\begin{document}

\title{Многомерные латинские битрейды \thanks{Работа выполнена при
поддержке Российского фонда фундаментальных исследований (проекты
11-01-00997, 10-01-00616) и ФЦП <Научные и научно-педагогические
кадры инновационной России> на 2009-2013 гг. (гос. контракт №
02.740.11.0362)}}
\author{В.\,Н.~Потапов}

\maketitle

\begin{center}
\textit{Институт математики им. С.\,Л.\,Соболева СО РАН,\\
Новосибирский государственный университет, Новосибирск }
\end{center}

\begin{abstract}

 Подмножество $k$-значного $n$-мерного гиперкуба
называется унитрейдом (объединённым битрейдом) если мощности  его
пересечений с одномерными гранями гиперкуба принимают только два
значения $0$ и $2$. Унитрейд называется двудольным (гамильтоновым),
если соответствующий ему подграф гиперкуба является двудольным
(гамильтоновым). Пара долей двудольного унитрейда называется
$n$-мерным латинским битрейдом. Для троичного $n$-мерного гиперкуба
определено число различных унитрейдов и получена экспоненциальная
нижняя оценка числа неэвивалентных латинских  битрейдов. Перечислены
все возможные $n$-мерные латинские битрейды мощности меньшей чем
$2^{n+1}$.

Подмножество $k$-значного $n$-мерного гиперкуба называется {
$t$-кратным МДР-кодом}, если оно пересекается с каждой одномерной
гранью гиперкуба ровно по $t$ элементам. Симметрическая разность
двух однократных МДР-кодов является  двудольным  унитрейдом. Каждая
из компонент соответствующего латинского битрейда является
свитчинговой компонентой одного из этих МДР-кодов. В статье
исследованы вопросы о мощностях компонент МДР-кодов и возможности
получения латинских битрейдов заданной мощности из МДР-кодов. Кроме
того, доказано, что любой МДР-код вкладывается в гамильтонов
двукратный МДР-код.

{\bf Ключевые слова}:  МДР-код, латинский битрейд, унитрейд,
компонента.

\end{abstract}

\section*{Введение}

Пусть $Q_k=\{0,1,\dots,k-1\}$. Обозначим через $Q_k^n$ множество
упорядоченных $k$-ичных наборов (вершин)
 длины $n$.
{\it Расстоянием Хэмминга} $d({x},{y})$ между
 вершинами $x,y\in Q_k^n $
  называется число позиций, в которых наборы ${x}$ и ${y}$
  различаются. Через  $\Gamma Q_k^n$ обозначим граф минимальных
расстояний метрического пространства $(Q_k^n, d)$. {\it Гранью
размерности} $k$ называется подмножество гиперкуба $Q_k^n$,
состоящее из вершин с одинаковыми фиксированными значениями
некоторых $n-k$ координат. В частности, {\it одномерная грань
направления} $i$, проходящая через вершину $(a_1,\dots,a_n)\in
Q_k^n$, определяется как множество
$\{(a_1,\dots,a_{i-1},x,a_{i+1},\dots,a_n) \ | \ x\in Q_k\}$.

  {\it {\rm МДР}-кодом} {\it с расстоянием $d$} называется
множество $M\subset Q_k^n$, пересекающееся с каждой $(d-1)$-мерной
гранью ровно по одному элементу.  Множество $W$ называется {\it
$t$-кратным МДР-кодом}, если оно пересекается с каждой одномерной
гранью ровно по $t$ элементам. Понятие кратного МДР-кода с
расстоянием $d$ совпадает с понятием корреляционо-иммунной функции
порядка $n-d+1$ в $Q_k^n$. Нетрудно видеть, что множество $M\subset
Q_k^n$ является МДР-кодом c расстоянием $d$ тогда и только тогда,
когда оно имеет мощность $k^{n-d+1}$ и расстояние между любыми двумя
различными элементами из $M$ не менее $d$. $t$-Кратный МДР-код
называется {\it расщепляемым}, если он является объединением $t$
однократных МДР-кодов. Кратные нерасщепляемые МДР-коды с расстоянием
$2$ рассматриваются в \cite{KP04}.
  В дальнейшем
в статье рассматриваются МДР-коды (в том числе кратные) только с
расстоянием $2$.

Граф называется {\it гамильтоновым}, если он содержит простой цикл,
проходящий через все вершины графа. Для однократного МДР-кода $M$
можно определить граф минимальных расстояний $\Gamma_2  M$, в
котором соединены ребром вершины, находящиеся на расстоянии $2$. Из
определения МДР-кода нетрудно получить, что граф $\Gamma_2  M$
всегда является гамильтоновым. Рассмотрим двукратный МДР-код
$D\subset Q_k^n$, через $\Gamma D$ обозначим граф минимальных
расстояний множества $D$, т.\,е. соединим рёбрами вершины,
находящиеся на расстоянии $1$. Доказано (\S\,3), что любой МДР-код
вкладывается в гамильтонов двукратный МДР-код. Заметим, что
двукратный МДР-код $D$ является расщепляемым тогда и только тогда,
когда граф $\Gamma  D$ является двудольным.

 Множество $B\subset Q_k^n$ будем
называть  {\it унитрейдом}\footnote{\,В \cite{KrMDS} использовался
термин 2-код.}, если мощности его пересечений с одномерными гранями
принимают только два значения $0$ и $2$. Унитрейд $B\subset Q_k^n$
будем называть {\it двудольным} унитрейдом, если подграф $\Gamma B$
графа $\Gamma Q_k^n$, порождённый множеством вершин $B$ является
двудольным. Ясно, что симметрическая разность двух МДР-кодов
является двудольным унитрейдом.

Поскольку в каждой одномерной грани МДР-кода содержится только одна
точка, можно считать, что МДР-код в $Q_k^n$ неявно задаёт некоторую
функцию от $n-1$ переменной. Таблица значений этой функции является
латинским ($n-1$)-мерным кубом порядка $k$ (латинским квадратом при
$n=3$), а МДР-код  можно рассматривать как график функции. Латинским
битрейдом называют пару частичных латинских квадратов, объединение
графиков которых является двудольным унитрейдом. Мы распространяем
название "латинский битрейд"\ на многомерный случай. Кроме того, в
тех случаях когда двудольный унитрейд однозначно разделяется на доли
(состоит из одной компоненты связности) будем называть латинским
битрейдом не только пару долей, но и сам унитрейд.

Двудольный унитрейд (латинский битрейд) $B$ будем называть {\it
полученным} из МДР-кода $M_1$ (возможно кратного), если найдётся
такой МДР-код $M_2$ (той же кратности), что $B=M_1\triangle M_2$. В
этом случае множество $B\cap M_1$ будем называть {\it компонентой}
МДР-кода $M_1$ в соответствии с общим представлением о свитчинговой
компоненте кода как о подможестве кода, замена которого на
равномощное сохраняет кодовое расстояние. Свитчинговым компонентам
МДР-кодов $M_1$ и $M_2$ соответствуют доли компонент связности графа
$\Gamma B$, они же компоненты латинского битрейда, соответствующего
унитрейду $B$.

В статье рассматриваются вопросы о числе унитрейдов и латинских
битрейдов в гиперкубе $Q_k^n$, о возможности получения латинских
битрейдов из МДР-кодов и о мощностях компонент МДР-кодов. В \S\,2
для троичного $n$-мерного гиперкуба определено число различных
унитрейдов ($2^{2^n}$) и получена асимптотическая (при
$n\rightarrow\infty$) нижняя оценка\footnote{\, Равенство
$f(n)=\Omega(g(n))$ означает, что $f(n)\geq cg(n)$ для некоторого
$c>0$ при $n\rightarrow\infty$.} $e^{\Omega(\sqrt{n})}$ числа
неэвивалентных латинских битрейдов. Нетривиальные верхние оценки
числа латинских битрейдов неизвестны. При $k\geq 4$ число и даже
асимптотика двойного логарифма числа унитрейдов и латинских
битрейдов также остаются неизвестными. Связь вопроса о числе
латинских битрейдов в гиперкубе $Q_k^n$
  с другими комбинаторными задачами
обсуждается в \cite{KP11}.

В \S\,1 перечислены все возможные $n$-мерные унитрейды мощности
меньшей чем $2^{n+1}$ и доказано,  что для каждого $s\in
\{0,\dots,n-1\}$ имеется единственный с точностью до эквивалентности
$n$-мерный латинский битрейд мощности $2^{n+1}-2^{s+1}$.  В \S\,3
доказана возможность получения этих латинских битрейдов  из
двукратных МДР-кодов и тем самым частично определён спектр мощностей
компонент двукратных МДР-кодов.

В  \cite{KP04}, \cite{Ph} имеется конструкция, позволяющая получать
из $t$-кратных МДР-кодов $t$-крат\-ные двоичные совершенные коды.
Вопросы связанные с мощностями компонент совершенных двоичных кодов
и гамильтоновости графов этих кодов рассматриваются в \cite{P12} и
\cite{Rm}. Полученные в данной работе результаты относительно
мощностей компонент МДР-кодов и гамильтоновости графов МДР-кодов
могут быть использованы при исследовании свойств других классов
кодов.

\section{Мощности унитрейдов}

\bpro\label{bitrpro0}

 a) Декартово произведение двух  унитрейдов является унитрейдом.

 b) Декартово произведение двух  двудольных унитрейдов является двудольным
унитрейдом.
 \epro

Доказательство пункта (a) непосредственно вытекает из определений,
при доказательстве пункта (b) пользуемся тем, что декартово
произведение двудольных графов является двудольным графом.

 Следующие два утверждения  доказаны в \cite{KrMDS} при $k=4$,
 однако, их можно доказать точно так же
 при произвольном $k\geq 2$.

\bpro\label{bitrpro1} Пусть $B\subset Q_k^n$ --- унитрейд, тогда
$|B|\geq 2^n$. \epro

\bpro\label{bitrpro11} Пусть $B\subset Q_k^n$ --- унитрейд.
Следующие условия эквивалентны:

a)  $|B|=2^n$,

b) для любой $m$-мерной грани мощность её пересечения с $B$
равняется $0$ или $2^m$,

c) унитрейд  $B$ пересекается только с двумя гипергранями каждого
направления,

d) подграф графа $\Gamma B$  изоморфен булеву кубу $\Gamma Q_2^n$.
\epro

Из предложения \ref{bitrpro1} непосредственно вытекает, что в
$Q_2^n$ существует единственный унитрейд, совпадающий со всем
множеством $Q_2^n$.

Рассмотрим вопрос о спектре мощностей  унитрейдов.

\bpro\label{bitrpro2} Пусть $B\subset Q_k^n$ --- унитрейд и
$2^{n+1}>|B|\geq 2^n$. Тогда $|B|=2^{n+1}- 2^{s+1} $, где $s\in
\{0,\dots,n-1\}$. \epro \proofr Будем доказывать утверждение методом
индукции по $n$. При $n=1$ утверждение очевидно, предположим оно
верно при $n-1$. Если унитрейд $B\subset Q_k^n$ содержится в
объединении двух гиперграней каждого направления, то   $|B|=2^n$ по
предложению \ref{bitrpro11}.

Если унитрейд $B$ пересекается с четырьмя гипергранями одного
направления, то по предложению \ref{bitrpro11} его мощность больше
или равна $2^{n+1}$.
 Пусть унитрейд $B$ пересекается с тремя  гипергранями
одного направления. Если пересечение хотя бы с одной из них имеет
мощность большую либо равную $2^{n}$, то по предложению
\ref{bitrpro1} имеем $|B|\geq 2^{n+1}$. В противном случае по
предположению индукции имеем $|B|=3\cdot2^n -
2^{s_1}-2^{s_2}-2^{s_3}$. Поскольку неравенство
$2^{s_1}+2^{s_2}+2^{s_3}>2^n$ выполнено только когда как минимум два
из трёх $s_i$ равняются $n-1$,  имеем $|B|=2^{n+1} - 2^{s}$.
\proofend

Как видно из доказательства предложения \ref{bitrpro2} унитрейд
$B\subset Q_k^n$ мощности меньшей чем $2^{n+1}$ пересекается не
более чем с тремя гипергранями любого направления, следовательно,
может быть вложен в троичный гиперкуб.

\bpro\label{bitrpro00}

 a) Симметрическая разность  двух  унитрейдов в гиперкубе $Q^n_3$ является
 унитрейдом.

b) Если симметрическая разность двух двудольных унитрейдов в
гиперкубе $Q^n_k$ является унитрейдом, а их пересечение   порождает
связный подграф в $Q^n_k$, то симметрическая разность является
латинским битрейдом. \epro

Доказательство пункта (a) непосредственно вытекает из определений,
при доказательстве пункта (b) пользуемся тем, что связный двудольный
граф однозначно разделяется на доли.

\bpro\label{bitrpro3} Для любого $s\in \{0,\dots,n-1\}$ существует
единственный (с точностью до эквивалентности)
  унитрейд $B_s\subset Q_3^n$ такой, что $|B_s|=2^{n+1}- 2^{s+1} $.
  Унитрейды $B_s$ являются латинскими битрейдами.
  \epro
\proofr  Из предложения \ref{bitrpro00} следует, что множество
\\ $B_s=( \{0,1\}^{n-s}\bigtriangleup\{1,2\}^{n-s})\times\{0,1\}^s$
является латинским битрейдом. Очевидно $|B|=2^s(2^{n-s}+2^{n-s}-2)$.

 Перейдём к доказательству
единственности. Рассмотрим унитрейд $B\subset Q^n_3$, $|B|<2^{n+1}$.
Аналогично доказательству предложения \ref{bitrpro2} получаем, что
пересечения множества $B$ с гипергранями некоторого направления
имеют мощности $2^{n-1}$, $2^{n-1}$, $2^n-2^s$. По предположению
\ref{bitrpro11} два пересечения эквивалентны булевым кубам, тогда из
определения унитрейда следует, что третье пересечение является
симметрической разностью  двух булевых кубов. \proofend

\section{Число унитрейдов в $Q_3^n$}

Рассмотрим множество функций $g:Q_3^n\rightarrow\{0,1\}$ как
векторное пространство $\mathbb{V}(n)$ над полем $GF(2)$.
Характеристические функции унитрейдов образуют в $\mathbb{V}(n)$
линейное подпространство $\mathcal{B}(n)$, в котором в качестве
базиса можно выбрать характеристические функции $\chi^B$ гиперкубов
вида $B=\{\alpha_1,2\}\times\dots\times\{\alpha_n,2\}$, где
$\alpha_i\in \{0,1\}$. Набор коэффициентов в разложении функции $g$
по базису $\{\chi^B\}$ является булевой функцией от набора
$(\alpha_1,\dots,\alpha_n)$. Далее в явной форме определим
преобразование, которое каждой булевой функции  ставит в
соответствие элемент из $\mathcal{B}(n)$.

  Определим частичный порядок на $Q_3$:
   $0<2$, $1<2$, а символы $0$ и $1$ несравнимы.
   Пусть $(x_1,\dots,x_n), (y_1,\dots,y_n)\in Q_3^n$.
   Введём обозначение $(x_1,\dots,x_n)\leq (y_1,\dots,y_n)$, если для любого
    $i\in \{1,\dots,n\}$ верно, что $x_i<y_i$ или
  $x_i=y_i$.
Отметим, что множество $\{x\in \{0,1\}^n\ |\ x\leq y\}$ является
гранью булева $n$-мерного гиперкуба размерности $wt(y)$ равной числу
 символов $2$ в наборе $y$. Более того, множество граней
находится во взаимно однозначном соответствии с множеством наборов
$y\in Q_3^n$.

  Пусть $f$ --- некоторая булева функция.
  Определим функцию $U[f]:Q_3^n\rightarrow \{0,1\}$
  равенством
$
 U[f](y)=
  \bigoplus\limits_{x\leq y}f(x) $.
Заметим, что в  булева функция $U[f]|_{\{0,2\}^n}$ является
преобразованием Мёбиуса от функции $f$.

\bpro\label{bitrpro31} a) Пусть $A\subseteq \{0,1\}^n$, тогда
$U[\chi^A]\in \mathcal{B}(n)$.

b) Пусть $g\in \mathcal{B}(n)$, тогда  $U[g|_{\{0,1\}^n}] =g$. \epro
  \proofr
a) Пусть $f=\chi^A$. Из определения преобразования $U$ имеем
равенство
$$U[f](a_1,\dots,a_{i-1},2,a_{i+1},\dots,a_n)= $$ $$=
U[f](a_1,\dots,a_{i-1},0,a_{i+1},\dots,a_n)\oplus
U[f](a_1,\dots,a_{i-1},1,a_{i+1},\dots,a_n)$$ для любых $a_j\in
Q_3$. Следовательно функция $U[f]$ имеет чётное число единиц в
каждой одномерной грани гиперкуба $Q_3^n$ и $U[f]\in
\mathcal{B}(n)$.

b) Равенство $U[g|_{\{0,1\}^n}](y) =g(y)$ для любого $y\in Q^n_3$
нетрудно  показать методом индукции по числу символов $2$ в наборе
$y$.
  \proofend

Из предложения \ref{bitrpro31} вытекает, что подпространство
$\mathcal{B}(n)$ имеет размерность $2^n$. Следовательно, справедливо
следующее

\bpro\label{bitrpro0} В гиперкубе $Q_3^n$ имеется ровно $2^{2^n}$
различных унитрейдов.\epro

Рассмотрим множество функций $g:Q_3^n\rightarrow R$ как векторное
пространство $\mathbb{V}_R(n)$ над полем $R$ вещественных чисел.
Рассмотрим линейное подпространство $\mathcal{B}_R(n)$, состоящее из
функций, сумма значений которых по любой одномерной грани равняется
$0$. Пусть $B\subset Q_3^n$
--- двудольный унитрейд. Определим функцию $h_B:Q^n_3 \rightarrow
\{-1,0,1\}$, которая принимает значение $1$ на первой доле
унитрейда, $-1$ --- на второй и $0$ --- в остальных вершинах
гиперкуба. Ясно, что $h_B\in \mathcal{B}_R(n)$. Определим оператор
$U_R$ для вещественнозначных функций аналогично оператору $U$. Пусть
$f:Q^n_2\rightarrow R$, тогда \begin{equation}\label{e1}
 U_R[f](y)=
 (-1)^{wt(y)} \sum\limits_{x\leq y}f(x).
 \end{equation}

Следующее предложение аналогично предложению \ref{bitrpro31}.

\bpro\label{bitrpro32} a) Для любой функции $f:Q_2^n\rightarrow R$
имеем $U_R[f]\in \mathcal{B}_R(n)$.

b) Пусть $g\in \mathcal{B}_R(n)$, тогда  $U_R[g|_{\{0,1\}^n}] =g$.
\epro

Пусть множество $\mathcal{F}(n)$ состоит из функций $f:Q^n_2
\rightarrow \{-1,0,1\}$, сумма значений которых в любой грани равна
одному из трёх чисел $-1,0,1$. {\it Счётчиком чётности} называется
булева функция $\delta(x_1,x_2,\dots,x_n)=\bigoplus_{i=1}^nx_i$.
Нетрудно видеть, что $(-1)^{\delta(x)}\in \mathcal{F}(n)$.

 Справедливо следующее

\bpro\label{bitrpro22} a) Пусть $B$ --- двудольный унитрейд, тогда
$h_B|_{\{0,1\}^n}\in \mathcal{F}(n)$.

b) Если $f\in \mathcal{F}(n)$, то $U_R[f]=h_B$ для некоторого
двудольного унитрейда $B\subset Q_3^n$. \epro \proofr

a) По предложению \ref{bitrpro32} (b) имеем $U_R[|h_B|_{\{0,1\}^n}]
=h_B$. Тогда  для любой грани $\{x\in \{0,1\}^n | x\leq y\}$ имеем
$\sum\limits_{x\leq y}h_B|_{\{0,1\}^n}(x)=(-1)^{wt(y)}h_B(y)\in
\{-1,0,1\}$.

b) По предложению \ref{bitrpro32} (a) имеем $U_R[f]\in
\mathcal{B}_R(n)$. Из условия $f\in \mathcal{F}(n)$ и определения
оператора $U_R$ следует, что $U_R[f](Q_3^n)\subseteq \{-1,0,1\}$.
Тогда в каждой одномерной грани функция $U_R[f]$ либо трижды
принимает значение $0$, либо по одному разу принимает значения
$-1,0,1$.
 \proofend

Рассмотрим некоторые конструкции функций из $\mathcal{F}(n)$.

\bpro\label{bitrpro23} a) Пусть $f\in \mathcal{F}(n)$, тогда $f\cdot
\chi^\gamma\in \mathcal{F}(n)$ для любой грани $\gamma$.

b) Пусть $\gamma_1,\gamma_2$ --- грани в $Q_2^n$ и
$\gamma_1\cap\gamma_2\neq \varnothing$. Определим функцию $f$
равенством
$$
f(x_1,\dots,x_n,x_{n+1}) =\begin{cases}
\chi^{\gamma_1}(-1)^{\delta(x_1,\dots,x_n)}   & \mbox{при } \
x_{n+1}=0 \cr
 \chi^{\gamma_2}(-1)^{\delta(x_1,\dots,x_n)\oplus1} & \mbox{при }\ x_{n+1}=1.
\end{cases}
$$
 Тогда
$f\in \mathcal{F}(n+1)$.

\epro

\bpro\label{bitrpro24}
 Пусть $f\in \mathcal{F}(n) $, $g\in \mathcal{F}(m) $ и
$F(x,y)=f(x)g(y)$. Тогда $F\in \mathcal{F}(n+m) $. \epro

Доказательства предложений \ref{bitrpro23} и \ref{bitrpro24}
нетрудно получить непосредственной проверкой.

Заметим, что все унитрейды в $Q^n_3$ состоят из одной компоненты
связности, поэтому двудольные унитрейды являются латинскими
битрейдами.

Далее получим нижнюю оценку числа неэквивалентных латинских
битрейдов. Две определённые на гиперкубе функции называются
эквивалентными, если они переходят друг в друга посредством
изометрий гиперкуба.
 Рассмотрим  гиперкуб $Q_2^n$ как векторное пространство над полем $GF(2)$.
Будем называть {\it носителем } вектора $x\in Q_2^n$ множество
позиций, на которых в векторе $x$ находятся единицы. Рассмотрим
набор векторов $z^1,\dots, z^k$ с попарно непересекающимися
носителями. Пусть $V\subset Q_2^n$ --- подпространство натянутое на
вектора $z^1,\dots, z^k$, $V=\{\bigoplus\alpha_iz^i \ | \alpha\in
Q_2^k\}$. Пусть $f:Q_2^k\rightarrow \{-1,0,1\}$. Определим функцию
$G_V[f]:Q_2^n\rightarrow \{-1,0,1\}$ равенствами

$$
G_V[f](x) =\begin{cases}f(\alpha)  & \mbox{при } \
x=\bigoplus\alpha_iz^i \cr
 0 & \mbox{при }\ x\not\in V .
\end{cases}
$$

\begin{theorem}
a) Если $f\in \mathcal{F}(k)$, то $G_V[f]\in \mathcal{F}(n)$.

b) Множество $\mathcal{F}(n)$ содержит не менее $e^{\Omega(\sqrt{
n})}$ неэквивалентных функций.

c) В гиперкубе $Q_3^n$ имеется не менее $e^{\Omega(\sqrt{ n})}$
неэквивалентных латинских битрейдов.
\end{theorem}
\proofr

a)  Поскольку носители векторов $z^1,\dots, z^k$ с попарно не
пересекаются, сумма значений функции $G_V[f]$ по  грани гиперкуба
$Q^n_2$ совпадает с суммой значений функции $f$ по некоторой грани
гиперкуба $Q^k_2$. Тогда $G_V[f](Q^n_3)\subseteq \{-1,0,1\}$.

b) Известно (см., например, \cite{Ed}), что имеется
$e^{\Omega(\sqrt{ n})}$ различных разбиений числа $n$ на целые
неотрицательные слагаемые. Каждое такое разбиение порождает
некоторый набор векторов $z^1,\dots, z^k$ с попарно
непересекающимися носителями и соответствующее набору
подпространство $V$. Функции $G_V[\delta]$ неэквивалентны для
различных разбиений, поскольку расстояния Хэмминга между базисными
векторами не меняются при изометриях гиперкуба.

c) По предложению \ref{bitrpro22} (b) каждая функция
$U_R[G_V[\delta]]$ порождает латинский битрейд $B$. Рассмотрим
случай, когда все носители векторов $z^1,\dots, z^k$ имеют мощность
не менее трёх и сумма этих мощностей  равна $n$. Покажем, что
пересечения битрейда $B$ c  гипергранями вида $\gamma=\{x\in Q^n_3 |
x_i=2\}$ обладают некоторым инвариантным относительно изометрии
свойством (*), которым не обладают  пересечения битрейда $B$ c
гипергранями вида $\{x\in Q^n_3 | x_i=0\}$ и $\{x\in Q^n_3 |
x_i=1\}$.  Тогда для любой изометрии $\varphi$ по битрейду
$\varphi(B)$ можно определить множество $\varphi(B\cap\{0,1\}^n)$.
Следовательно, неэквивалентность битрейдов $U_R[G_V[\delta]]$ при
различных $V$ является следствием пункта b).

(*) Пусть $\gamma=\{x\in Q^n_3 | x_i=a\}$, $a\neq 2$ тогда и только
тогда, когда найдётся такая гипергрань $\gamma'=\{x\in Q^n_3 |
x_j=b\}$, $b\in Q_3$, $j\neq i$, что пересечение
$\gamma\cap\gamma'\cap B$ пусто.

Без ограничения общности полагаем, что $i=1$ и носитель вектора
$z^1$ содержит $1$ и $2$. Тогда
$$\{x\in Q^n_3 | x_1=0\}\cap \{x\in Q^n_3 |
x_2=1\}\cap B=\varnothing, $$

$$\{x\in Q^n_3 | x_1=1\}\cap \{x\in Q^n_3 |
x_2=0\}\cap B=\varnothing. $$

В то время как из формулы (\ref{e1}) имеем

$$y=(2,0,\dots,0) \in \{x\in Q^n_3 | x_1=2\}\cap \{x\in Q^n_3 |
x_j=0\}\cap B, $$

$$y'=(2,0,\dots,0,2,0\dots) \in \{x\in Q^n_3 | x_1=2\}\cap \{x\in Q^n_3 |
x_j=2\}\cap B, $$

поскольку $\{\overline{0}\}=$ $$=\{x\in \{0,1\}^n  | x\leq y\}\cap
\{x\in \{0,1\}^n  | G_V[\delta](x)\neq 0\}= \{x\in \{0,1\}^n  |
x\leq y'\}\cap \{x\in \{0,1\}^n | G_V[\delta](x)\neq 0\}.$$

Кроме того,
$$y''=(2,0,\dots,0,1,\dots,1,0\dots) \in \{x\in Q^n_3 | x_1=2\}\cap \{x\in Q^n_3 |
x_j=1\}\cap B, $$

поскольку\\  $\{x\in \{0,1\}^n \ | x\leq y''\}\cap \{x\in \{0,1\}^n
\ | \ G_V[\delta](x)\neq 0\}
=\{(0,\dots,0,1,\dots,1,0\dots)\}=\{z^k\}$, где носитель вектора
$z^k$ содержит позицию $j$.
 \proofend

Отметим, что латинский битрейд $B_s$ (см. предложение
\ref{bitrpro3}) можно задать также функцией
$h_{B_s}=U_R[G_V[\delta]]$, где базис подпространства $V$ состоит из
одного вектора $z$ с носителем мощности $n-s$.

\section{МДР-коды}

Методом индукции по размерности $n$ нетрудно доказать, что любой
МДР-код в гиперкубе $Q^n_3$ продолжается до МДР-кода в гиперкубе
$Q^{n+1}_3$ ровно двумя способами, причём при каждом $n\geq 1$ все
МДР-коды в гиперкубе $Q^n_3$ эквивалентны (см. \cite{LM}, Exercise
13.15).
 Характеризация всех МДР-кодов в  гиперкубе  $Q_4^n$
имеется в \cite{KP09}.  Как было указано выше, симметрическая
разность $M_1\triangle M_2$ двух МДР-кодов $M_1$ и $M_2$ является
латинским битрейдом, причём каждая из долей $M_1\cap(M_1\triangle
M_2)$ и $M_2\cap(M_1\triangle M_2)$ является компонентой
соответствующего МДР-кода. Таким образом, в предложении
\ref{bitrpro2} получены ограничения на возможные мощности компонент
МДР-кодов. Отметим (см., например, \cite{KP11}), что в гиперкубах
$Q_k^n$ при $k\geq 4$ имеются МДР-коды с  компонентами мощности
кратной $2^{n-1}$ и соответстующие им унитрейды мощности кратной
$2^{n}$. Более того, в \cite{KrMDS} показано, что любой унитрейд,
содержащийся в двукратном МДР-коде в $Q^n_4$, имеет мощность кратную
$2^n$. Ниже рассмотрен вопрос о возможности получения из МДР-кодов
(в том числе кратных)  унитрейдов других мощностей.

 Для произвольной функции $f:Q_k^n\rightarrow Q_k$
обозначим её график через ${\cal M}\langle f\rangle = \{(x,f( x)):
x\in Q_k^n\}$. Если множество ${\cal M}\langle f\rangle $ является
МДР-кодом, то функция $f$ называется $n$-арной квазигруппой порядка
$k$. Соответственно, частичной $n$-арной квазигруппой порядка $k$
называется функция, график которой пересекается с любой одномерной
гранью не более чем по одной вершине.

В \cite{Cruse} доказано, что любая частичная $n$-арная квазигруппа
конечного порядка есть сужение  $n$-арной квазигруппы некоторого
большего порядка. Отсюда следует

\bpro\label{bitrpro4} Любой  латинский битрейд $B\subset Q_k^n$
может быть получен из некоторого МДР-кода $M\subset Q_m^n$, где
$m\geq k$. \epro

\bpro\label{bitrpro5}
 Латинский битрейд $B\subset Q_k^n$ такой, что $2^{n+1}>|B|> 2^n$,
 не может быть получен из МДР-кода  $M\subset Q_k^n$ при $k\in\{3,4\}$. \epro
\proofr
 Поскольку все МДР-коды в $Q_3^n$ эквивалентны при любом
$n$,
 достаточно рассмотреть
произвольный МДР-код в $Q_3^n$, например, $M=\{x\in Q_3^n\ |\
x_1+\dots+x_n=0\, {\rm mod }\,3\}$. Нетрудно видеть, что любой
полученный из $M$ унитрейд содержит весь МДР-код $M$ и имеет
мощность $2\cdot3^{n-1}$.

Пусть $k=4$.  Рассмотрим трёхмерные латинские битрейды. Полным
перебором нетрудно убедиться, что если полученный из МДР-кода
$M\subseteq Q_4^3$ битрейд $B$  имеет пересечение мощности не менее
$6$ с некоторой двумерной гранью, то найдётся ещё одна двумерная
грань параллельная первой, пересечение с которой также имеет
мощность не менее $6$. Кроме того, полным перебором можно проверить,
что в этом случае $|B|\geq 16 =2^4$.

Рассмотрим латинский битрейд $B'$, полученный из МДР-кода в $Q^n_4$.
Применяя метод индукции, получаем, что если битрейд $B'\subseteq
Q_4^n$ имеет пересечение мощности не менее $6$ с какой-либо
двумерной гранью, то $|B'|\geq 2^{n+1}$. Из предложения
\ref{bitrpro11} следует, что если пересечение латинского битрейда с
любой двумерной гранью имеет мощность меньше $6$, т.\,е. равную $0$
или $4$, то граф порождённый битрейдом  $B'$ является булевым
$n$-кубом.\proofend

Рассмотрим возможность получения латинских битрейдов малой мощности
из двукратных МДР-кодов.

\bpro\label{bitrpro6} Любой $n$-мерный латинский  битрейд мощности
меньше $2^{n+1}$ может быть получен из двукратных МДР-кодов в
$Q_4^{n}$. \epro \proofr Вначале покажем, что битрейд
$B_0^k=\{0,1\}^{k}\bigtriangleup\{1,2\}^{k}$ может быть получен из
некоторого двукратного МДР-кода в $Q_4^{k}$. Пусть функция
$g:Q^n_k\rightarrow\{0,1\}$ определена равенством
$g(a_1,\dots,a_k)=\sum\limits_{i=1}^k a_i \ {\rm mod
  }\,2$. Нетрудно видеть, что функции $g$, $g_1=g\oplus \chi^{\{1,2\}^k}$ и
   $g_2=g\oplus\chi^{\{0,1\}^k}$
  являются характеристическими функциями некоторых двукратных МДР-кодов $M_0$, $M_1$
 и $M_2$. Поскольку $g_1\oplus g_2 = \chi^{\{1,2\}^k} \oplus
 \chi^{\{0,1\}^k}$, имеем $B_0^k=M_1\bigtriangleup M_2$.

Рассмотрим булевозначные функции

\begin{equation}\label{be2}
 f_1(a_1,\dots,a_k,y)=\begin{cases} g_1(a_1,\dots,a_k)  & \mbox{при }\ y=0,\cr
 g_2(a_1,\dots,a_k)  & \mbox{при }\ y=1, \\
g_1(a_1,\dots,a_k)\oplus1  & \mbox{при }\ y=2,\\
 g_2(a_1,\dots,a_k)\oplus1  & \mbox{при }\ y=3;
\end{cases}
\end{equation}

и

\begin{equation}\label{be3}
 f_2(a_1,\dots,a_k,y)=\begin{cases} g_2(a_1,\dots,a_k)  & \mbox{при }\ y=0,\cr
 g_1(a_1,\dots,a_k)  & \mbox{при }\ y=1, \\
g_1(a_1,\dots,a_k)\oplus1  & \mbox{при }\ y=2,\\
 g_2(a_1,\dots,a_k)\oplus1  & \mbox{при }\ y=3.
\end{cases}
\end{equation}

Легко видеть, что $f_1$ и $f_2$ являются характеристическими
функциями двукратных МДР-кодов, причём $f_1\oplus f_2 =
\chi^{B^k_0\times\{0,1\}}$.

Подставляя в формулы (\ref{be2}--\ref{be3}) функции $f_1$ и $f_2$
вместо $g_1$ и $g_2$, обнаружим, что латинский битрейд
${B^k_0\times\{0,1\}^2}$ также может быть получен из некоторого
двукратного МДР-кода в $ Q_4^{k+2}$ и т.\,д. Для окончания
доказательства остаётся сослаться на предложение \ref{bitrpro3}.
  \proofend

Приведём  таблицы функций $g$, $g_1$, $g_2$  при $n=2$:

$   \left(
 \begin{array}{cccc}
0 & 1  & 0 & 1\\
1 & 0 & 1 & 0\\
0 & 1  & 0 & 1\\
1 & 0  & 1 & 0\\
 \end{array} \right),
$ $   \left(
 \begin{array}{cccc}
0 & 1  & 0 & 1\\
1 & 1 & 0 & 0\\
0 & 0  & 1 & 1\\
1 & 0  & 1 & 0\\
 \end{array} \right),
$$   \left(
 \begin{array}{cccc}
1 & 0  & 0 & 1\\
0 & 1 & 1 & 0\\
0& 1  & 0 & 1\\
1 & 0  & 1 & 0\\
 \end{array} \right).
$

 Можно показать, что построенные в предложении
\ref{bitrpro6} двукратные МДР-коды  являются нерасщепляемыми при
$n\geq 3$.

Далее исследуем вопрос о гамильтоновости графов минимальных
расстояний  двукратных МДР-кодов.

\bpro\label{bitrpro16} Для любого МДР-кода $M\subset Q_k^n$ найдётся
двукратный МДР-код $D$ такой, что  $M\subset D$ и граф минимальных
расстояний $\Gamma D$ является гамильтоновым.
 \epro
\proofr Определим МДР-код $M'\subset Q_k^n$ равенством $$
M'=\{(x_1+1{\rm mod
  }\,k,x_2,\dots,x_n)\ |\ (x_1,\dots,x_{n})\in M\}. $$
Рассмотрим двукратный МДР-код $D=M\cup M'$. Докажем  по индукции,
что граф $\Gamma D$ является гамильтоновым, причём для любого ребра
найдётся проходящий через него гамильтонов цикл. При $n=2$ граф
$\Gamma D$ представляет собой гамильтоновов цикл, состоящий из
чередующихся рёбер двух направлений.

Пусть предположение индукции доказано для при $n-1$. Рассмотрим
произвольное ребро в графе $\Gamma D$, например, ребро $v$
соединяющее точки $(0,\dots,0,c)$ и $(0,\dots,0,c')$.  Граф $\Gamma
D'$, где $D'=\{(x_1,x_{n}) \ |\ (x_1,0,\dots,0,x_{n})\in D \}$,
является простым циклом $H'$ и содержит ребро $v$. По предположению
индукции граф $\Gamma D_a$, где $D_a=\{(x_1,\dots,x_{n-1}) \ |\
(x_1,\dots,x_{n-1},a)\in D \}$, содержит гамильтонов цикл $H_a$ для
любого $a\in Q_k$, причём гамильтонов цикл $H_a$ можно выбрать
проходящим через ребро, принадлежащее циклу $H'$.

Если два простых цикла пересекаются ровно по одному ребру, то удалив
из  объединения двух циклов это ребро мы получим простой цикл. Таким
образом, граф $H=  H' \triangle \bigcup_a H_a$ является простым
циклом. Кроме того, цикл $H$ проходит через все вершины графа
$\Gamma D$ и содержит ребро $v$.
 \proofend

  \begin{center}

\end{center}
\end{document}